\documentclass[12pt]{amsart}

\usepackage{amsmath, amscd, amssymb}

\newcommand{\bC}{{\mathbb C}}

\newcommand{\bP}{{\mathbb P}}

\newcommand{\bR}{{\mathbb R}}
\newcommand{\bZ}{{\mathbb Z}}
\newcommand{\bn}{{\bf n}}
\newcommand{\bp}{{\bf p}} 
 \newcommand{\bw}{{\bf w}}

 \newcommand{\cP}{{\mathcal P}}

\newcommand{\Mbar}{\overline{\mathcal M}}

\newcommand{\pd}{\partial}
\newcommand{\half}{\frac{1}{2}}

\newcommand{\vmu}{{\vec{\mu}} }

\newcommand{\cor}[1]{\langle {#1} \rangle}
\DeclareMathOperator{\res}{Res}
\DeclareMathOperator{\Aut}{Aut}

\newtheorem{theorem}{Theorem}[section]
\newtheorem{Theorem}{Theorem}

\theoremstyle{remark}
\newtheorem{remark}{Remark}[section]

\theoremstyle{definition}
\newtheorem{example}{Example}[section]

\newcommand{\bea}{\begin{eqnarray}}
\newcommand{\eea}{\end{eqnarray}}
\newcommand{\ben}{\begin{eqnarray*}}
\newcommand{\een}{\end{eqnarray*}}
\newcommand{\be}{\begin{equation}}
\newcommand{\ee}{\end{equation}}

\begin{document}

\title
{Local Mirror Symmetry for the Topological Vertex}

\author{Jian Zhou}
\address{Department of Mathematical Sciences\\Tsinghua University\\Beijing, 100084, China}
\email{jzhou@math.tsinghua.edu.cn}

\begin{abstract}
For three-partition triple Hodge integrals related to the topological vertex,
we derive Eynard-Orantin type recursion relations from the cut-and-join equation.
This establishes a version of local mirror symmetry for the local $\bC^3$ geometry with
three $D$-branes,
as proposed by Mari\~no \cite{Mar} and Bouchard-Klemm-Mari\~no-Pasquetti \cite{BKMP}.
\end{abstract}

\maketitle

\section{Introduction}

This is a sequel to an earlier paper \cite{Zho} in which we study  a version of local  mirror symmetry
for one-legged framed topological vertex.
In this paper we generalize the results to the general topological vertex \cite{AKMV}.
The motivations of this work and also that of \cite{Zho}
can be found in the Introduction of that paper,
so we will only briefly explain them here.
Physically the topological vertex is a formalism
that computes partition functions in the local A-theory of certain noncompact Calabi-Yau $3$-folds with D-branes,
based on duality with Chern-Simons theory developed in a series of work.
A mathematical theory of the topological vertex \cite{LLLZ} has been developed,
based on some earlier work  on Hodge integrals and localizations
on relative moduli spaces.
Recently a new formalism  for the local B-theory on the mirror of toric Calabi-Yau threefolds has been proposed
in \cite{Mar, BKMP},
inspired by the recursion procedure of \cite{EO} discovered first in the context of matrix models.
It is then very interesting to verify local mirror symmetry in arbitrary genera
using this new formalism of the local B-theory,
and this has been done in many cases in \cite{Mar} and \cite{BKMP}.
For the simplest case which is the one-legged framed topological vertex,
Bouchard and Mari\~no \cite{Bou-Mar} made a conjecture
based on the proposed new formalism of the B-theory in \cite{Mar, BKMP}.
They also made a similar conjecture for Hurwitz numbers,
which has recently  been proved by Borot-Eynard-Safnuk-Mulase \cite{BEMS}
and Eynard-Safnuk-Mulase \cite{EMS} by two different methods.
More recently,
two slightly different proofs for the Bouchard-Mari\~no Conjecture for the one-legged framed
topological vertex have appeared \cite{Che, Zho}.
They both use ideas from  \cite{EMS}.
In this paper we will make a generalization of \cite{Zho} to study the case of general
topological vertex.
The following is our main result:

\begin{Theorem}
The three-partition triple Hodge integrals related to the topological vertex
satisfy some Eynard-Orantin type recursion relations,
as proposed by Mari\~no \cite{Mar} and Bouchard-Klemm-Mari\~no-Pasquetti \cite{BKMP}.
\end{Theorem}

Indeed, such recursion relations can be derived from the cut-and-join equation
satisfied by these Hodge integrals by the same method as in \cite{Zho}.

The rest of this paper is arranged as follows.
In \S 2 we elaborate on some constructions in the study of
the local mirror symmetries of $\bC^3$.
For example,
we show that some Lagrangian submanifolds correspond to complex submanifolds in the mirror manifold,
but some do not.
In \S 3 we derive some recursion relations for the three-partition triple Hodge integrals related to the topological
vertex using a cut-and-join equation they satisfy.
Finally in \S 4, we reformulate the recursion relations in terms of Eynard-Orantin type
recursions following the proposal of Mari\~no \cite{Mar} and Bouchard-Klemm-Mari\~no-Pasquetti \cite{BKMP}.

\vspace{.1in}
{\em Acknowledgements}.
This research is partially supported by two NSFC grants (10425101 and 10631050)
and a 973 project grant NKBRPC (2006cB805905).

\section{The Mirror Geometry of $\bC^3$}

For general cases of local mirror constructions,
see \cite{HV, HIV, AV, BKMP}.
In this section we will go through some details of the  constructions for $\bC^3$.

\subsection{T-duality on flat $3$-tori}

Let $V$ be a Euclidean space of dimension $3$,
and let $\Gamma \subset V$ be a lattice.
Then $V/\Gamma$ is a flat $3$-torus.
Let $V^*$ be the dual space of $V$, endowed with the dual Euclidean metric.
Let $\Gamma^* \subset V^*$ be the normalized dual lattice,
i.e.,
\be
\Gamma^* = \{ \varphi \in V^* \,|\, \frac{1}{4\pi^2} \langle\varphi, v \rangle \in \bZ, \;\forall v \in \Gamma\}.
\ee
The torus $V^*/\Gamma^*$ will be referred to as the {\em T-dual}of $V/\Gamma$.
Now because $(V^*)^*\cong V$ and $(\Gamma^*)^* \cong \Gamma$,
the T-dual of $V^*/\Gamma^*$ is $V/\Gamma$.

\begin{example} \label{exm:Lattice}
Let $V$ be $\bR^3$ with the standard metric,
$\Gamma$ be the lattice generated by $(2\pi r_1, 0, 0)$,
$(0, 2\pi r_2, 0)$, and $(0,0, 2\pi r_3)$.
Then $V/\Gamma$ is a Riemannian product of three copies of circles,
with radii $r_1$, $r_2$ and $\frac{r_3}{2\pi}
$, respectively.
The dual space $V^*$ is still $\bR^3$ with standard metric,
$\Gamma^*$ is now the lattice generated by $(2\pi/ r_1, 0, 0)$,
$(0, 2\pi /r_2, 0)$, and $(0,0, 2\pi /r_3)$,
and $V^*/\Gamma^*$ is a Riemannian product of three copies of circles,
with radii $1/r_1$, $1/r_2$ and $1/r_3$, respectively.
\end{example}

Let $\Gamma_1 \subset \Gamma$ be a sublattice,
which spans a linear subspace $V_1 \subset V$,
then $V_1/\Gamma_1$ is a subtorus of $V/\Gamma$.
One can get a family of $2$-tori by translating the subspace $V_1$ by vectors in $V$:
Given any $v \in V_1$, each affine subspace $v + V_1$ determines a subtorus in $V/\Gamma$.
Let $\Gamma_1^\perp \subset \Gamma^*$ be the sublattice of $\Gamma^*$,
orthogonal to $\Gamma_1$,
i.e.,
\be
\Gamma_1^\perp = \{ \varphi \in \Gamma^* \;|\; \varphi(v) = 0, \; \forall v \in \Gamma_1\}.
\ee
It is easy to see that the linear subspace spanned by $\Gamma_1^\perp$ in $V^*$
is the subspace $V_1^\perp$,
which is the orthogonal complement of $V_1$ in $V^*$.
The substorus $V_1^\perp/\Gamma_1^\perp$ will be referred as the {\em T-dual}
of $V_1/\Gamma_1 \subset V/\Gamma$.
We can also obtain from it a family of subtori by translating $V_1^\perp$ by vectors in $V^*$.

Note $(\Gamma_1^\perp)^\perp = \Gamma_1$ and $(V_1^\perp)^\perp = V_1$,
the T-dual of $V_1^\perp/\Gamma_1^\perp \subset V^*/\Gamma^*$ is $V_1/\Gamma_1$.

\begin{example}
Let $V$ and $\Gamma$ be as in Example \ref{exm:Lattice}.
Let $\Gamma_1$ be the sublattice generated by $(2\pi r_1, 4\pi r_2, 0)$.
Then $\Gamma_1^\perp$ is the sublattice generated by
$(-4\pi /r_2, 2\pi /r_1, 0)$ and $(0,0,2\pi/r_3)$.
In this case $V_1/\Gamma_1$ is a circle,
and $V_1^\perp/\Gamma_1^\perp$ is a $2$-torus.
\end{example}

\subsection{A degenerate torus fibration of $\bC^3$}
Consider the following natural torus action on $U(1)^3 \times \bC^3 \to \bC^3$:
\be \label{eqn:Action}
(e^{i\alpha_1}, e^{i\alpha_2}, e^{i\alpha_3}) \cdot
(z_1, z_2,z_3) = (e^{i\alpha_1} z_1, e^{i\alpha_2} z_2, e^{i\alpha_3} z_3).
\ee
This is a Hamiltonian action with moment map  given by:
\be
\mu: \bC^3 \to \bR_{\geq 0}^3, \quad \mu(z_1, z_2, z_3) = \half(|z_1|^2,|z_2|^2, |z_3|^2),
\ee
where $\bR_{\geq 0}$ is the half line $\{x\in \bR \;|\; x \geq 0\}$.
The inverse image of $\mu$ is generically a flat $3$-torus,
but along the boundary of $\bR_{\geq 0}^3$,
it may degenerate to a $2$-torus, a circle or a point.
Let $\bR_+ =\{x \in \bR\;|\;x > 0\}$.
We will focus on $\mu: (\bC^*)^3 \to \bR_+^3$,
where $\bC^*$ is the set of nonzero complex numbers.
This is a fibration with flat $3$-tori as fibers.
By adding the degenerate fibers corresponding to the points on $\pd \bR_{\geq 0}^3$,
one gets a partial compactification:
The union of the degenerate fibers form a divisor with normal crossing in $\bC^3$:
\ben
\{(z_1,z_2, 0)\;|\;z_1,z_2\in \bC\} \cup \{(z_1,0, z_3)\;|\;z_1,z_3\in \bC\} \cup
\{(0,z_2, z_3)\;|\;z_2,z_3\in \bC\}.
\een

On $(\bC^*)^3$ one can use the polar coordinates on each copies of $\bC^*$,
then the Euclidean metric becomes
\be
g = \sum_{j=1}^3 (dr_j^2 + r_j^2d\theta_j^2).
\ee

\subsection{T-duality of the torus fibration}
\label{sec:Mirror}

Note we have a fibration of flat $3$-tori,
so we can apply the T-duality to each fiber.
More precisely,
one can define a new metric on $(\bC^*)^3$:
\be
g^\vee = \sum_{j=1}^3 (dr_j^2 + \frac{1}{r_j^2}d\theta_j^2).
\ee
In this metric,
the $3$-tori in the fibers becomes their dual tori.
This metric is no longer Hermitian with respect to the original complex structure.
For the metric $\hat{g}$ to be a Hermitian metric,
the dual almost complex structure is taken to be:
\begin{align} \label{eqn:JDual}
J^\vee(\frac{\pd}{\pd r_j}) & = r_j \frac{\pd}{\pd \theta_j}, &
J^\vee (\frac{\pd}{\pd \theta_j}) & = - \frac{1}{r_j} \frac{\pd}{\pd r_j}.
\end{align}
On differential $1$-forms,
the almost complex structure acts by:
\begin{align}
dr_j & \mapsto - \frac{1}{r_j} d\theta_j, & d\theta_j & \mapsto r_j dr_j.
\end{align}
Note the action of the original almost complex structure acts as follows on tangent vectors and $1$-forms:
\begin{align}
\frac{\pd}{\pd r_j} & \mapsto \frac{1}{r_j} \frac{\pd}{\pd \theta_j}, &
\frac{\pd}{\pd \theta_j} & \mapsto - r_i \frac{\pd}{\pd r_j}, \\
dr_j & \mapsto - r_j d\theta_j, & d\theta_j & \mapsto \frac{1}{r_j} dr_j.
\end{align}
In the dual almost complex structure,
$\half r_j^2 + \sqrt{-1} \theta_ji$ are new complex local coordinates,
because $r_j dr_j + \sqrt{-1} d\theta_j$ are now of type $(1,0)$.
Because $\theta_j$ is multiple valued,
we take
\be \label{Def:Y}
y_j = \exp (-(\half r_j^2 + \sqrt{-1} \theta_j)).
\ee
Denote  by $X^\vee$ the space $(\bC^*)^3$ with the dual complex coordinates $\{y_1, y_2, y_3\}$.
Note because $|y_j|<1$,
$X^\vee$ is no longer a complex $3$-torus $(\bC^*)^3$,
but instead a product of three copies of punctured open unit disc
$$\{y \in \bC\;|\; 0< |y| < 1\}.$$

We will denote this dual complex space by $X^\vee$.
In the dual complex coordinates the metric $\hat{g}$ can be written as:
\be
\hat{g} = - \sum_{j=1}^3 \frac{|dy_j|^2}{|y_j|^2 \cdot \ln |y_j|^2}.
\ee
Note this is a complete metric, but no longer Ricci-flat.
The dual symplectic form is given by
\be
\omega^\vee = -\sum_{j=1}^3 \frac{1}{r_j} dr_j \wedge d \theta_j
= -\frac{\sqrt{-1}}{2} \sum_{j=1}^3 \frac{dy_j \wedge d \bar{y}_j }{|y_j|^2 \cdot \ln |y_j|^2}.
\ee

Now the $3$-torus action (\ref{eqn:Action}) is given in polar coordinates by
\begin{align}
r_j & \mapsto r_j, & \theta_j & \mapsto \theta_j + \alpha_j.
\end{align}
Therefore,
it induces the following action on the mirror manifold:
\be
y_j \mapsto e^{i\alpha_j} y_j.
\ee
This is again a Hamiltonian action with moment map
$\mu^\vee: X^\vee \to \bR^3$ given by:
\be
\mu^\vee(y_1, y_2, y_3)
= (\ln \ln |y_1|^2, \ln \ln |y_2|^2, \ln \ln |y_3|^2).
\ee

\subsection{Lagrangian submanifolds of Aganagic-Vafa type}
We recall some Lagrangian submanifolds of $\bC^3$ constructed by Aganagic and Vafa \cite{AV}.
Let
\bea
&& m_1 \mu_1 + m_2 \mu_2 +  m_3 \mu_3 = \alpha_1,  \label{eqn:AV1} \\
&& n_1 \mu_1 + n_2 \mu_2 +  n_3 \mu_3 = \alpha_2 \label{eqn:AV2}
\eea
be two planes in $\bR_+^2$ that intersects along a ray $l$ in $\bR_+^3$.
We assume that $m_i, n_i$ are integers.
This ray can also be described by:
\be \label{eqn:Ray}
(\mu_1, \mu_2, \mu_3) = (a_1, a_2, a_3) + t_1 (k_1, k_2, k_3),
\ee
for some integers $k_1, k_2, k_3$ and nonnegative integers $a_1, a_2, a_3, t_2$.
For $b_1, b_2, b_3 \in \bR$,
\begin{align}
\theta_i = m_i s + n_i t + b_i, \qquad i=1,2,3, \qquad s, t\in \bR \label{eqn:AV3}
\end{align}
determines a subtorus $K$ of the $3$-torus.
It is clear that $l \times K$ is a Lagrangian submanifold in $\bC^3$.

\begin{example} \label{exm:AV}
Consider the following complex conjugate involution on $(\bC^*)^3$:
\be
(z_1, z_2, z_3) \mapsto (\frac{a^2}{\bar{z}_1}, \bar{z}_3, \bar{z}_2),
\ee
for some $a > 0$.
The fixed point set is a Larangian submanifold given by
\be
\{(a e^{i\theta_1}, r_2 e^{i\theta_2}, r_2 e^{-i\theta_2}) \; |\; \theta_1, \theta_2 \in \bR, \,r_2 \in \bR_+\}.
\ee
This is a Lagragian submanifold  of Aganagic-Vafa type which corresponds to the following ray in $\bR^3_+$:
$$(\half a^2, t, t), \quad t > 0.$$
A Lagrangian submanifold will be said to be in phase I, phase II or phase III if it
corresponds to a ray of the form $(a, t, t)$, $(t,a,t)$ or $(t,t,a)$ for fixed $a > 0$, respectively.
\end{example}

\begin{example} \label{exm:KL}
One can also consider the following complex conjugate involution \cite{Kat-Liu}:
\be
(z_1, z_2, z_3) \mapsto (\frac{1}{\bar{z}_1}, \bar{z}_1\bar{z}_3, \bar{z}_1\bar{z}_2).
\ee
The fixed point set is a Larangian submanifold given by
\be
\{(e^{i\theta_1}, r_2 e^{i\theta_2}, r_2e^{-i(\theta_1+\theta_2)}) \; |\; \theta_1, \theta_2 \in \bR, \,r_2 \in \bR_+\}.
\ee
This is not of Aganagic-Vafa type.
\end{example}

\subsection{Mirror $B$-branes of Lagrangian submanifold of Aganagic-Vafa type}
One can apply the T-duality to each of the fibers in a Lagrangian submanifold
of Aganagic-Vafa type,
because they are subtori of the $3$-tori.
It was claimed in \cite{AV} that this will yield a complex submanifold.
We now verify this claim by the discussion of the dual almost complex structure
presented in \S \ref{sec:Mirror}.
Take  a dual subtorus $K$ of dimension $2$ in $T^3$ as in (\ref{eqn:AV3}).
Its dual can be taken a circle given in the theta coordinates by
\bea
&& m_1 \theta_1 + m_2 \theta_2 +  m_3 \theta_3 = \beta_1,  \label{eqn:AV1Dual} \\
&& n_1 \theta_1 + n_2 \theta_2 +  n_3 \theta_3 = \beta_2, \label{eqn:AV2Dual}
\eea
or equivalently
\be \label{eqn:AV3Dual}
(\theta_1, \theta_2, \theta_3) = (c_1, c_2, c_3) + t_2(k_1, k_2, k_3).
\ee
Combined with (\ref{eqn:AV3}),
one sees that the tangent space is given by the linear space of the following two vectors:
\begin{align*}
\sum_{j=1}^3 k_j \frac{\pd}{\pd r_j}, && \sum_{j=1}^3 k_j \frac{\pd}{\pd \theta_j}
\end{align*}
This is $J^\vee$-invariant,
therefore, we have proved the following

\begin{theorem}
The mirror dual $l \times K^\vee$ of a Lagrangian submanifold $l \times K$ of Aganagic type
is a complex submanifold of $X^\vee$.
\end{theorem}

\begin{remark}
It is straightforward to extend this result to the case of $\bC^n$ for arbitrary $n$ and substorus of
any codimension.
\end{remark}

In our case, $l \times K^\vee$ can be explicitly given in the $y_j$-coordinates.
By (\ref{eqn:AV3}) and (\ref{eqn:AV3Dual}):
\be
y_j = e^{-(a_j + ic_j)} e^{- k_j (t_1 + i t_2)},
\ee
where $t_1, t_2 > 0$,
or by (\ref{eqn:AV1}), (\ref{eqn:AV2}), (\ref{eqn:AV1Dual}), (\ref{eqn:AV2Dual}),
we have
\bea
&& y_1^{m_1} y_2^{m_2} y_3^{m_3} = e^{- (\alpha_1 + i \beta_1)}, \\
&& y_1^{n_1} y_2^{n_2} y_3^{n_3} = e^{- (\alpha_2 + i \beta_2)}.
\eea

\begin{example}
The mirror dual of the Lagrangian submanifold in Example \ref{exm:AV} is given by:
\begin{align}
y_1 & = e^{-(a^2/2+ c_1 i)}, & y_2 & = e^{-c_2i} e^{-(t_1+it_2)}, & y_3 & = e^{-c_3 i} e^{-(t_1+it_2)}.
\end{align}
\end{example}

\begin{example}
The mirror dual of the Lagrangian submanifold in Example \ref{exm:KL} is given by
\begin{align*}
r_1 & = 1, & r_2 & = r_3 = r, & \theta_1 & =\theta_2=\theta_2 = \theta.
\end{align*}
Therefore, its tangent space is
spanned by the following vectors;
\begin{align*}
\frac{\pd}{\pd r_2} + \frac{\pd}{\pd r_3}, && \sum_{j=1}^3  \frac{\pd}{\pd \theta_j}
\end{align*}
This is not $J^\vee$-invariant,
and so $l \times K^\vee$ is not a complex submanifold of $X^\vee$.
\end{example}

\subsection{The mirror curve of $\bC^3$}

The mirror geometry of $\bC^3$ according to the construction of Hori-Iqbal-Vafa \cite{HIV}
is not the mirror manifold $X^\vee$ above,
but instead the subspace of $\bC^2 \times (\bC^*)^3$ defined by the following equation
\be
uv = y_1+y_2+y_3
\ee
where $u,v \in \bC$, $y_1, y_2 \in \bC^*$,
modulo the following action by $\bC^*$:
\be
t \cdot (u, v, y_1, y_2, y_3) \mapsto (tu, v, ty_1, ty_2, ty_3).
\ee
Denote this space by $\hat{X}$,
and by $[u,v,y_1, y_2, y_3]$ the equivalence class of a point $(u, v, y_1, y_2, y_3)$ in this space.
There is a natural projection map $\pi: \hat{X} \to \bP^2$ defined by
\be
\pi([u,v, y_1, y_2, y_3]) = [y_1:y_2:y_3] \in \bP^2.
\ee
The base space of this projection is the projectivized $X^\vee$:
\be
\{[y_1:y_2:y_3] \in \bP^2\;|\; y_1, y_2, y_3 \in (\bC^*)^2\},
\ee
and the fiber $\pi^{-1}[y_1, y_2, y_3])$ is the space
\be
\{ (\tilde{u}, v) \in \bC^2\;|\; \tilde{u} v = \frac{y_1}{y_3} + \frac{y_2}{y_3} + 1 \}.
\ee
This is copy of $\bC^*$ when the RHS does not vanish,
but it becomes a copy of normal crossing singular set $\{(\tilde{u},v)\in \bC^2\;|\;\tilde{u} v = 0\}$
when
\be
y_1 + y_2 + y_3 = 0.
\ee
One can embed $(\bC^*)^3$ in this space as follows:
\be
(y_1, y_2, y_3) \mapsto (w_1, w_2, w_3, u, v)
= (y_1, y_2, 1, y_3, \frac{y_1+y_2+1}{y_3}),
\ee
The curve
\be
\{[y_1:y_2:y_3] \in \bP^2\;|\; y_1, y_2, y_3 \in \bC^*, \; y_1 + y_2 + y_3 = 0\}
\ee
is called the {\em mirror curve} of  $\bC^3$.
This is a copy of $\bC^*$,
or equivalently, $\bP^1$ with three points $0, -1, \infty$ removed.
If we take $z=\frac{y_1}{y_3}$ as coordinate on the mirror curve,
$z = 0, 1, \infty$ correspond to $[y_1:y_2: y_3] = [0: -1: 1]$, $[-1: 1; 0]$ and
$[1:-1:0]$ respectively.
There are different ways to realize the mirror curve as a plane curve.
If one takes
\begin{align}
\tilde{X} & = \frac{y_1}{y_3}, & \tilde{Y} & = \frac{y_2}{y_3},
\end{align}
then one gets
\be
\tilde{X} + \tilde{Y} + 1 = 0;
\ee
if one takes
\begin{align}
X' & = \frac{y_3}{y_1}, & Y' & = \frac{y_2}{y_1},
\end{align}
then one gets
\be
X'+Y'+1=0
\ee
and
\begin{align}
X' & = \tilde{X}^{-1}, & Y' & = \tilde{X}^{-1} \tilde{Y}.
\end{align}
One can also take
\begin{align}
X & = \tilde{x}\tilde{y}^{a}, & Y & = \tilde{y}
\end{align}
for any integer $a$, then one gets:
\be
X + Y^a  +  Y^{a+1}= 0.
\ee
Let $X=-(-1)^a x$ and $Y = -Y$,
then one gets:
\be \label{eqn:FramedMirror}
x = y^a - y^{a+1}.
\ee
This is the equation for the framed mirror curve.

\section{Recursion Relations from Cut-and-Join equation for the Topological Vertex}
\label{sec:CJ}

In this section we will first recall the definition of the three-partition triple Hodge integrals
related to the topological vertex.
Then we will derive some recursions relations
using the cut-and-join equation.

\subsection{Partitions}
Recall that a partition $\mu$ of a nonnegative integer $d$ is
a sequence of positive integers
$\mu=(\mu_1\geq \mu_2 \geq \cdots \geq \mu_n>0)$ such that
$$ d=\mu_1+\ldots+\mu_n.
$$
We write $|\mu|=d$, and  $l(\mu)=n$.
The following numbers associated with a partition $\mu$ are often used:
\bea
&& |\Aut(\mu)|=\prod_{j}m_j(\mu)!,\quad\text{where}\quad
m_j(\mu)=\#\{i:\mu_i=j\}, \\
&& z_\mu= |\Aut(\mu)| \cdot \prod_{i=1}^n \mu_i, \\
&& \kappa_\mu=\sum_{i=1}^{l(\mu)}\mu_i(\mu_i-2i+1),
\eea
For convenience we also consider the partition of $0$ and write it as $\emptyset$.
We use the following conventions:
\begin{align}
|\Aut(\emptyset)| & = 1, & z_{\emptyset} & = 1, & \kappa_{\emptyset} & = 0.
\end{align}

Let $\cP$ denote the set of partitions, and let
$$
\cP_+=\cP-\{\emptyset\},\ \
\cP^2_+=\cP^2-\{(\emptyset,\emptyset)\},\ \
\cP^3_+=\cP^3-\{(\emptyset,\emptyset,\emptyset)\}.
$$
Given a triple of partitions $\vmu=(\mu^1, \mu^2,\mu^3) \in\cP^3$, we define
$$
l(\vmu)=\sum_{i=1}^3 l(\mu^i),\ \
\Aut(\vmu)=\prod_{i=1}^3\Aut(\mu^i).
$$

\subsection{Three-Partition Hodge Integrals related to the topological vertex}

Let $T^3$ act on $\bC^3$ by
\be
(e^{iw_1}, e^{iw_2}, e^{iw_3}) \cdot (z_1, z_2, z_3)
= (e^{iw_1} z_1, e^{iw_2} z_2, e^{iw_3} z_3).
\ee
Elements of the subgroup of $T^3$ that preserves the holomorphic volume form $dz_1 \wedge dz_2 \wedge dz_3$
satisfy the following condition:
\be
w_1+w_2+w_3=0.
\ee
We will focus on this subgroup and assume this condition throughout the rest of this section.
For simplicity of notations,
we will write:
\begin{equation}
\bw=(w_1, w_2, w_3), \quad w_4=w_1.
\end{equation}
We will let
\be \label{eqn:a}
a_i: = \frac{w_{i+1}}{w_i},
\ee
and write $a_1 = a$.
Then we have:
\begin{align} \label{eqn:a2a3}
\frac{w_3}{w_2} & =- \frac{a+1}{a}, & \frac{w_1}{w_3} & = - \frac{1}{a+1}.
\end{align}
For $\vmu=(\mu^1, \mu^2, \mu^3)\in \cP^3_+$, we let
$$
d^1_\vmu=0,\quad
d^2_\vmu= l(\mu^1),\quad
d^3_\vmu= l(\mu^1)+ l(\mu^2).
$$
We define the {\em three-partition triple Hodge integral}
to be
\begin{multline} \label{eqn:Gg}
G_{g;\vmu}(\bw) = \frac{(-\sqrt{-1})^{l(\vmu)} }{|\Aut(\vmu)|}
 \prod_{i=1}^3\prod_{j=1}^{l(\mu^i)}
\frac{\prod_{k=1}^{\mu^i_j-1}(\mu^i_j w_{i+1} + k w_i) }
     {(\mu^i_j-1)!w_i^{\mu^i_j-1} } \\
\int_{\Mbar_{g, l(\vmu)}
} \prod_{i=1}^3\frac{\Lambda_g^\vee(w_i)w_i^{l(\vmu)-1} }
{\prod_{j=1}^{l(\mu^i)}(w_i(w_i-\mu^i_j\psi_{d^i_\vmu+j}))},
\end{multline}
where $\Lambda_g^\vee(u)=u^g-\lambda_1 u^{g-1}+\cdots + (-1)^g\lambda_g$.
See \cite{LLLZ} for its relationship with the topological vertex.
From the definition,
we have the following cyclic symmetries:
\be \label{eqn:Cyclic}
G_{g; \mu^1, \mu^2, \mu^3}(w_1, w_2, w_3)
= G_{g; \mu^2, \mu^3, \mu^1}(w_2, w_3,w_1).
\ee
Note $\sqrt{-1}^{l(\vmu)} G_{g,\vmu}(\bw)$ is a rational function in
$w_1, w_2, w_3$ with rational coefficients,
and by a simple dimension counting is homogeneous of degree $0$.
We will write
\be
G_{g,\vmu}(a) : = G_{g,\vmu}(\bw) = G_{g, \mu}(1, a, -1-a).
\ee

Several exceptional cases that play important roles need special care.
Recall for $n \geq 3$,
the following identity holds:
\be
\int_{\Mbar_{0,n}} \frac{1}{(1-a_1\psi_1) \cdots (1-a_n\psi_n)}
= (a_1+ \cdots + a_n)^{n-3}.
\ee
This identity inspires the following useful conventions:
\bea
&& \int_{\Mbar_{0,2}} \frac{1}{(1-a_1\psi_1) (1-a_2\psi_2)}
= \frac{1}{a_1+ a_2}, \\
&& \int_{\Mbar_{0,2}} \frac{1}{1-a_1\psi_1}
= \frac{1}{a_1^2}.
\eea
By these conventions we have in the case of $l(\vmu) = 1$,
\ben
&& G_{0; (m), \emptyset, \emptyset}(a) = -\sqrt{-1} \frac{\prod_{j=1}^{m-1} (ma+j)}{(m-1)!} \frac{1}{m^2},
\een
and similar expressions for $G_{0; \emptyset, (m), \emptyset}(a)$
and $G_{0; \emptyset, \emptyset, (m)}(a)$ by changing $a$ to $a_2 = - \frac{a+1}{a}$
and $a_3 = - \frac{1}{a+1}$ respectively.
In the case of $l(\vmu) = 2$,
we have
\ben
&& G_{0; (m_1, m_2), \emptyset, \emptyset}(a) =  \frac{a(a+1)}{1+\delta_{m_1, m_2}}
\prod_{i=1}^2 \frac{\prod_{j=1}^{m_i-1} (m_ia+j)}{(m_i-1)!} \cdot \frac{1}{m_1+m_2},
\een
and similar expressions for $G_{0; \emptyset, (m_1,m_2), \emptyset}(a)$
and $G_{0; \emptyset, \emptyset, (m_1,m_2)}(a)$ by changing $a$ to $- \frac{a+1}{a}$ and $- \frac{1}{a+1}$ respectively;
furthermore,
\bea
&& G_{0; (m^1_1), (m^2_1),\emptyset}(a)
= -\frac{1}{a_3} \prod_{i=1,2} \frac{\prod_{j=1}^{m^i_1-1} (m^i_1 a_i +j)}{(m^i_1-1)!} \cdot \frac{1}{m^1_1 a+m^2_1}, \\
&& G_{0; \emptyset, (m^2_1), (m^3_1)}(a)
= -\frac{1}{a_1} \prod_{i=2,3} \frac{\prod_{j=1}^{m^i_1-1} (m^i_1 a_i +j)}{(m^i_1-1)!} \cdot \frac{1}{m^2_1 a_2 +m^3_1}, \\
&& G_{0; (m^1_1),\emptyset, (m^3_1)}(a)
= - \frac{1}{a_2} \prod_{i=3,1} \frac{\prod_{j=1}^{m^i_1-1} (m^i_1 a_i +j)}{(m^i_1-1)!} \cdot \frac{1}{m^3_1 a_3 +m^1_1}.
\eea

\subsection{The cut-and-join equation for three-partition Hodge integrals}

Let $p^i=(p^i_1,p^i_2,\ldots)$ be formal variables.
Given a partition $\mu$ we define $p^i_\mu=p^i_1\cdots p^i_{\ell(\mu)}$; (note
$p^i_\emptyset=1$).
We write
$$
\bp = (p^1, p^2, p^3), \qquad \bp_\vmu=p^1_{\mu^1}p^2_{\mu^2}p^3_{\mu^3}.
$$
Define the generating functions of the three-partition Hodge integrals by
\ben
&& G_\vmu(\lambda;a)=\sum_{g=0}^\infty
\lambda^{2g-2+\ell(\vmu)}G_{g,\vmu}(a), \\
&& G(\lambda;\bp;a)=\sum_{\vmu\in\cP^3_+} G_\vmu(\lambda;a)\bp_\vmu,
\een
where $\lambda$ is a formal variable.
By the results in \cite{LLLZ},
the following cut-and-join equation is satisfied by $G$:
\begin{multline} \label{eqn:CJ}
 \frac{\pd G}{\pd a}
= \frac{\sqrt{-1}\lambda}{2} \sum_{k=1}^3 \frac{\pd}{\pd a} (\frac{w_{k+1}}{w_k}) \\
\cdot \sum_{i, j \geq 1} \biggl(
ijp^k_{i+j} \frac{\pd^2 G}{\pd p^k_i\pd p^k_j}
+ ijp^k_{i+j}\frac{\pd G}{\pd p^k_i}\frac{\pd G}{\pd p^k_j}
+ (i+j)p^k_ip^k_j\frac{\pd G}{\pd p^k_{i+j}} \biggr).
\end{multline}
By (\ref{eqn:a}) and (\ref{eqn:a2a3}),
\begin{align}
\frac{\pd}{\pd a} (\frac{w_2}{w_1}) & = 1, &  \frac{\pd}{\pd a} (\frac{w_3}{w_2}) & = \frac{1}{a^2}, &
\frac{\pd}{\pd a} (\frac{w_1}{w_3}) & = \frac{1}{(a+1)^2}.
\end{align}

\subsection{The symmetrized generating function for three-partition Hodge integrals}

One can also define
\ben
&& G_g(\bp;a)
= \sum_{\vmu} G_{g; \vmu}(a) \bp_{\vmu}.
\een
Because $G_g(\bp;a)$ is a formal power series in $p^k_1, p^k_2, \dots, p^k_n, \dots$,
$k=1,2,3$,
for each $\bn=(n_1, n_2, n_3)$,
one can obtain from it a formal power series $\Phi_{g;\bn}(x^1_{[n_1]}; x^2_{[n_2]}; x^3_{[n_3]}; a)$
by applying the following symmetrization operator \cite{Gou-Jac-Vai, Che}:
$$p_{\mu^k}^k \mapsto (\sqrt{-1})^{n_k}
\delta_{l(\mu^k), n_k}\sum_{\sigma \in S_{n_k}}
(x^k_{\sigma(1)})^{\mu^k_1} \cdots (x^k_{\sigma(n_k)})^{\mu^k_{n_k}}.$$

\begin{multline}
G_{g;\vmu}(\bw) = \frac{(-\sqrt{-1})^{\ell(\vmu)} }{|\Aut(\vmu)|}
 \prod_{i=1}^3\prod_{j=1}^{\ell(\mu^i)}
\frac{\prod_{a=1}^{\mu^i_j-1}(\mu^i_j w_{i+1} + a w_i) }
     {(\mu^i_j-1)!w_i^{\mu^i_j-1} } \\
\int_{\Mbar_{g,\ell(\vmu)}
} \prod_{i=1}^3\frac{\Lambda_g^\vee(w_i)w_i^{\ell(\vmu)-1} }
{\prod_{j=1}^{\ell(\mu^i)}(w_i(w_i-\mu^i_j\psi_{d^i_\vmu+j}))}
\end{multline}

From the definition,
we have for $2g-2+n_1+n_2+n_3 > 0$,
\be \label{eqn:PhiInX}
\begin{split}
& \Phi_{g;\bn}(x^1_{[n_1]}; x^2_{[n_2]}; x^3_{[n_3]}; a) \\
= &  (-a(a+1))^{n_1+n_2+n_3-1} \sum_{b^i_j \geq 0}
\cor{\prod_{i=1}^3 \prod_{j=1}^{n_i} \tau_{b^i_j} \cdot T_g(a)}_g
\prod_{i=1}^3 \prod_{j=1}^{n_i} \phi_{b^i_j}(x^i_j; a_i) \\
& \cdot \prod_{j=1}^{n_2}  \frac{1}{a^{2+b^2_j}} \cdot \prod_{j=1}^{n_3}  \frac{1}{(-a-1)^{2+b^3_j}} ,
\end{split}
\ee
where
\ben
&& \cor{\tau_{b_1} \cdots \tau_{b_n} T_g(a)}_g
= \int_{\Mbar_{g,n}} \prod_{i=1}^n \psi_i^{b_i} \cdot
\Lambda_g^{\vee}(1)\Lambda_g^{\vee}(a)\Lambda_g^{\vee}(-1-a), \\
&& \phi_b(x;a) =  \sum_{m\geq 1}\frac{\prod_{j=1}^{m-1}(ma+j)}{(m-1)!} m^b x^m.
\een
For simplicity of notations we will set
\bea
&& \phi^1_b(x;a) = \phi_b(x;a), \\
&& \phi^2_b(x;a_2) = \frac{1}{a^{2+b^2}} \phi_b(x;a_2), \\
&& \phi^3_b(x;a_3) = \frac{1}{(-a-1)^{2+b}} \phi_b(x;a_3).
\eea

\subsection{Exceptional cases}
We now study several exceptional cases that will play a key role later.
First we have
\be \label{eqn:Phi01}
\Phi_{0;1,0,0}(x^1_1;a) = \sum_{m=1}^{\infty} \frac{\prod_{j=1}^{m-1} (m a +j)}{(m-1)!} m^{-2} (x_1^1)^m
= \phi_{-2}(x_1^1;a);
\ee
similarly,
\bea
&& \Phi_{0;0,1,0}(x^2_1;a) = \phi_{-2}(x_1^2;a_2), \\
&& \Phi_{0;0,0,1}(x^3_1;a) = \phi_{-2}(x_1^3;a_3).
\eea
So we have
\bea
&& x^1 \frac{\pd}{\pd x^1} \Phi_{0;1,0,0}(x_1^1;a)
= \phi_{-1}(x^1_1;a) = -\ln y(x^1_1;a),  \label{eqn:Phi01(1)} \\
&& x^2 \frac{\pd}{\pd x^2} \Phi_{0;0,1,0}(x_1^2;a)
= \phi_{-1}(x^2_1;a_2) = -\ln y(x^2_1;a_2), \label{eqn:Phi01(2)} \\
&& x^3 \frac{\pd}{\pd x^3} \Phi_{0;0,0,1}(x_1^3;a)
= \phi_{-1}(x^3_1;a_3) = -\ln y(x^3_1;a_3), \label{eqn:Phi01(3)}
\eea
where
\be
y(x;a) = 1 - \sum_{n=1}^{\infty} \frac{\prod_{j=0}^{n-2} (na+j)}{n!} x^n.
\ee
Secondly,
\ben \label{eqn:Phi02}
\Phi_{0;2,0,0}(x^1_1, x^1_2; a)
= -a(a+1) \sum_{m_1, m_2 \geq 1} \prod_{i=1}^2 \frac{\prod_{j=1}^{m_i-1} (m_ia+j)}{(m_i-1)!}
\cdot \frac{(x^1_1)^{m_1}(x_2^1)^{m_2}}{m_1+m_2}.
\een
This has been treated in \cite{Che, Zho}.
It is easy to see that:
\begin{multline} \label{eqn:EqnPhi(02)}
\big(x^1_1\frac{\pd}{\pd x^1_1} + x^1_2\frac{\pd}{\pd x^1_2}\big) \Phi_{0;2,0,0}(x^1_1, x^1_2;a) \\
= -a (a+1) \frac{t(x_1^1;a)-1}{a+1} \cdot \frac{t(x^1_2;a)-1}{a+1} \\
= - a(a+1) \frac{y_1-1}{(a+1) y_1 - a} \cdot \frac{y_2-1}{(a+1) y_2 - a}.
\end{multline}
One can verify that:
\begin{multline} \label{eqn:Phi02(1)}
\Phi_{0;2,0,0}(x^1_1, x^1_2; a)\\
=-\ln(\frac{y(x^1_2;a)-y(x^1_1;a)}{x^1_1-x^1_2})
+ \ln \frac{1-y(x^1_1;a)}{x^1_1}+\ln \frac{1-y(x^1_2;a)}{x^1_2},
\end{multline}
\begin{multline} \label{eqn:Phi02(1)Der1}
x^1_1 \frac{\pd}{\pd x^1_1} \Phi_{0;2,0,0}(x^1_1,x^1_2;a) \\
= -(\frac{1}{y(x^1_1;a)-y(x^1_2;a)} - \frac{1}{y(x^1_1;a)-1} )
x^1_1\frac{\pd y(x^1_1;a)}{\pd x^1_1}+\frac{x^1_2}{x^1_1-x^1_2}.
\end{multline}
\begin{multline} \label{eqn:Phi02(1)Der2}
\frac{\pd}{\pd x_1^1} \frac{\pd}{\pd x_2^1} \Phi_{0;2,0,0}(x^1_1,x^1_2;a) \\
=-\frac{1}{(y(x^1_1;a)-y(x^1_2;a))^2} \frac{\pd y(x^1_1;a)}{\pd x^1_1}
\frac{\pd y(x^1_2;a)}{\pd x^1_2} + \frac{1}{(x^1_1-x^1_2)^2}.
\end{multline}
One can get similar results for
$\Phi_{0;0,2,0}(x^2_1, x^2_2; a)$ and $\Phi_{0;0,0,2}(x^3_1, x^3_2; a)$
by changing the indices and changing $a$ to $a_2$ and $a_3$ respectively.
The case of $\Phi_{0;1,1,0}(x^1_1; x^2_1; a)$, $\Phi_{0;0,1,1}(x^2_1; x^3_1; a)$
 and $\Phi_{0;1,0,1}(x^1_1; x^3_1; a)$ can be treated in the same fashion.
First,
\ben \label{eqn:Phi12Series}
\Phi_{0;1,1,0}(x^1_1; x^2_1; a)
= - (a+1) \sum_{m^1, m^2 \geq 1} \prod_{i=1}^2 \frac{\prod_{j=1}^{m^i-1} (m^ia_i+j)}{(m^i-1)!}
\cdot \frac{(x^1_1)^{m^1}(x_1^2)^{m^2}}{m^1a+m^2}.
\een
It is easy to see that:
\begin{multline}
\big(a x^1_1\frac{\pd}{\pd x^1_1} + x_1^2\frac{\pd}{\pd x_1^2}\big) \Phi_{0;1,1,0}(x^1_1; x_1^2;a) \\
= - (a+1) \frac{t(x_1^1;a)-1}{a+1} \cdot \frac{t(x_1^2;a_2)-1}{a_2+1} \\
= - (a+1) \frac{y(x_1^1;a)-1}{(a+1) y(x_1^1;a) - a} \cdot \frac{y(x_1^2;a_2)-1}{(a_2+1) y(x_1^2;a_2) - a_2}.
\end{multline}
Write $y_1 = y(x_1^1;a)$ and $y_2 = y(x_1^2;a_2)$,
one then gets:
\be \label{eqn:EqnPhi(12)}
\big(a x^1_1\frac{\pd}{\pd x^1_1} + x_1^2\frac{\pd}{\pd x_1^2}\big) \Phi_{0;1,1,0}(x^1_1; x_1^2;a)
= \frac{a(y_1-1)(y_2-1)}{(y_1 - \frac{a}{a+1})(y_2 - (a+1))}.
\ee
One can verify that:
\be \label{eqn:Phi(12)}
\Phi_{0;1,1,0}(x_1^1; x_1^2; a)
= - \ln \frac{\tilde{y}(x^1_1; a) - y(x_1^2;a_2)}{\tilde{y}(x_1^1;a)-1}.
\ee
Indeed,
\begin{multline} \label{eqn:Phi(12)Eqn}
- (a x_1^1 \frac{\pd}{\pd x_1^1} + x_1^2 \frac{\pd}{\pd x_1^2})
\ln \frac{\tilde{y}(x^1_1; a) - y(x_1^2;a_2)}{\tilde{y}(x_1^1;a)-1} \\
= - \frac{a x_1^1 \frac{\pd \tilde{y}(x_1^1;a)}{\pd x_1^1} - x_1^2 \frac{\pd y(x_1^2; a_2)}{\pd x_1^2}  }{\tilde{y}(x^1_1; a) - y(x_1^2;a_2)}
+ \frac{ a x_1^1 \frac{\pd \tilde{y}(x_1^1;a)}{\pd x_1^1}}{\tilde{y}(x_1^1;a)-1}.
\end{multline}
Recall
\ben
x_1^2 \frac{\pd y(x_1^2; a_2)}{\pd x_1^2}
= \frac{y(x_1^2;a_2)(1- y(x_1^2;a_2))}{a_2 - (a_2+1) y(x_1^2;a_2)}
= - \frac{a y_2(1-y_2)}{(a+1) - y_2}.
\een
Now we use $x_1^1 = y(x_1^1;a)^{a}(1- y(x_1^1;a))$ to get:
\be
\tilde{y}(x_1^1;a)
= \frac{1}{x_1^1} y(x_1^1;a)^a = \frac{1}{1 - y(x_1^1; a)} = \frac{1}{1-y_1}.
\ee
It follows that:
\ben
&& x_1^1 \frac{\pd \tilde{y}(x_1^1;a)}{\pd x_1^1}
= x_1^1 \frac{\pd}{\pd x_1^1} ( \frac{1}{1 - y(x_1^1; a)} ) \\
& = &  \frac{1}{(1 - y(x_1^1; a) )^2} x_1^1\frac{\pd}{\pd x_1^1}  \tilde{y}(x_1^1;a) \\
& = & \frac{y(x_1^1; a)}{(1 - y(x_1^1;a)) \cdot (a - (a+1) y(x_1^1;a))}
 = \frac{y_1}{(1-y_1)(a-(a+1)y_1)}.
\een
Therefore,
the RHS of (\ref{eqn:Phi(12)Eqn}) is
\ben
&& - \frac{\frac{ay_1}{(1-y_1)(a-(a+1)y_1)} + \frac{ay_2(1-y_2)}{a+1-y_2}}{\frac{1}{1-y_1} - y_2}
+ \frac{\frac{ay_1}{(1-y_1)(a-(a+1)y_1)} }{\frac{1}{1-y_1} - 1} \\
& = & \frac{a(y_1-1)(y_2-1)}{(y_1-\frac{a}{a+1})(y_2-(a+1))}.
\een
This matches with the RHS of (\ref{eqn:EqnPhi(12)}).
From (\ref{eqn:Phi(12)}) one easily gets:
\be \label{eqn:Phi(12)Der1}
x^1_1 \frac{\pd}{\pd x_1^1}  \Phi_{0;1,1,0}(x_1^1; x_1^2; a)
= - \frac{x_1^1 \frac{\pd \tilde{y}(x_1^1;a)}{\pd x_1^1}   }{\tilde{y}(x^1_1; a) - y(x_1^2;a_2)}
+ \frac{x_1^1 \frac{\pd \tilde{y}(x_1^1;a)}{\pd x_1^1}}{\tilde{y}(x_1^1;a)-1},
\ee
\be \label{eqn:Phi(12)Der2}
\frac{\pd}{\pd x_1^1} \frac{\pd}{\pd x_1^2} \Phi_{0;1,1,0}(x_1^1; x_1^2; a)
= - \frac{\frac{\pd \tilde{y}(x^1_1;a)}{\pd x_1^1} \frac{\pd y(x_1^2;a_2)}
{\pd x^2_1}}{(\tilde{y}(x^1_1;a) - y(x_1^2;a_2))^2}.
\ee
By the cyclic symmetries one gets from (\ref{eqn:Phi(12)}):
\be \label{eqn:Phi(23)}
\Phi_{0;0,1,1}(x_1^2; x_1^3; a)
= - \ln \frac{\tilde{y}(x^2_1; a_2) - y(x_1^3;a_3)}{\tilde{y}(x_1^2;a_2)-1},
\ee
and
\be \label{eqn:Phi(31)}
\Phi_{0;1,0,1}(x_1^1; x_1^3; a)
= - \ln \frac{\tilde{y}(x^3_1; a_3) - y(x_1^1;a_1)}{\tilde{y}(x_1^3;a_3)-1}.
\ee
And so we have
\be \label{eqn:Phi(23)Der2}
\frac{\pd}{\pd x_1^2} \frac{\pd}{\pd x_1^3} \Phi_{0;0,1,1}(x_1^2; x_1^3; a)
= - \frac{\frac{\pd \tilde{y}(x^2_1;a_2)}{\pd x_1^2} \frac{\pd y(x_1^3;a_3)}{\pd x^3_1}}
{(\tilde{y}(x^2_1;a_2) - y(x_1^3;a_3))^2},
\ee
and
\be \label{eqn:Phi(31)Der1}
x^1_1 \frac{\pd}{\pd x^1_1} \Phi_{0;1,0,1}(x_1^1; x_1^3; a)
=  \frac{x^1_1 \frac{\pd y(x^1_1;a)}{\pd x^1_1}}{\tilde{y}(x^3_1; a_3) - y(x_1^1;a_1)},
\ee
\be \label{eqn:Phi(31)Der2}
\frac{\pd}{\pd x_1^1} \frac{\pd}{\pd x_1^3} \Phi_{0;1,0,1}(x_1^1; x_1^3; a)
= - \frac{\frac{\pd \tilde{y}(x^3_1;a_3)}{\pd x_1^3} \frac{\pd y(x_1^1;a)}{\pd x^1_1}}
{(\tilde{y}(x^3_1;a_3) - y(x_1^1;a))^2}.
\ee

\begin{remark}
Without the proposal in \cite{BKMP},
it will be very difficult for the author to find the explicit expressions (\ref{eqn:Phi(12)}),
(\ref{eqn:Phi(23)}), and (\ref{eqn:Phi(31)}).
\end{remark}

\subsection{The symmetrized cut-and-join equation for three-partition Hodge integrals}

Using the analysis in \cite{Gou-Jac-Vai},
one can obtain the symmetrized version of the cut-and-join equation (\ref{eqn:CJ}) as follows:
\be \label{eqn:CJinX}
\begin{split}
&  \frac{\pd}{\pd a} \Phi_{g; \bn}(x^1_{[n_1]}; x^2_{[n_2]}; x^3_{[n_3]}; a) \\
=&  \frac{1}{2}\sum_{k=1}^{n_1} z_1 \frac{\pd}{\pd z_1}
z_2 \frac{\pd}{\pd z_2}\Phi_{g-1; n_1+1, n_2, n_3}(z_1, z_2,x^1_{[n_1]_k}; x^2_{[n_2]}; x^3_{[n_3]}; a)|_{z_1, z_2=x_k^1}\\
+ & \frac{1}{2} \sum_{k=1}^{n_1}
\sum_{\substack{g_1+g_2 = g \\
A^1 \coprod B^1 = [n_1]_k \\A^2 \coprod B^2 = [n_2] \\ A^3 \coprod B^3 = [n_3]}}
x^1_k \frac{\pd}{\pd x^1_i} \Phi_{g_1; |A^1|+1, |A^2|, |A^3|}(x^1_k, x^1_{A^1};x^2_{A^2}; x^3_{A^3}; a) \\
& \cdot
x^1_k \frac{\pd}{\pd x^1_i} \Phi_{g_2; |B^1|+1, |B^2|, |B^3|}(x^1_k, x^1_{B^1}; x^2_{B^2}; x^3_{B^3}; a) \\
- & \sum_{k=1}^{n_1} \sum_{j \in [n_1]_k}  \frac{x^1_j}{x^1_k-x^1_j} \cdot
x^1_k \frac{\pd}{\pd x^1_k} \Phi_{g; n_1-1, n_2, n_3}(x^1_{[n_1]_j}; x^2_{[n_2]}; x^3_{[n_3]};  a) \\
+ & \cdots.
\end{split}
\ee
Here we have omitted the terms that correspond to cut-and-join operation in $x^2_i$ and $x^3_i$.

\subsection{Symmetrized cut-and-join equation in the $v$-coordinates}
As in \cite{EMS, Che, Zho},
introduce the $v$-coordinates by
\be
x_j^i = e^{-(v_j^i)^2/2}.
\ee
Now for $2g-2+n_1+n_2+n_3> 0$,
(\ref{eqn:PhiInX}) becomes
\be \label{eqn:PhiInV}
\begin{split}
& \Phi_{g;\bn}(x^1_{[n_1]}; x^2_{[n_2]}; x^3_{[n_3]}; a) \\
= &  (-a(a+1))^{n_1+n_2+n_3-1} \sum_{b^i_j \geq 0}
\cor{\prod_{i=1}^3 \prod_{j=1}^{n_i} \tau_{b^i_j} \cdot T_g(a)}_g
\prod_{i=1}^3 \prod_{j=1}^{n_i} \xi^i_{b^i_j}(v^i_j; a_i)
\end{split}
\ee
where
\bea
&& \xi^1_b(v;a) = \xi_b(v;a), \\
&& \xi^2_b(v;a) = \frac{1}{a^{2+b}} \xi_b(v;a_2), \\
&& \xi^3_b(v;a) = \frac{1}{(-a-1)^{2+b}} \xi_b(v;a_3).
\eea
The first term on the right-hand side is now:
\ben
&& \frac{1}{2}\sum_{k=1}^{n_1} z_1 \frac{\pd}{\pd z_1}
z_2 \frac{\pd}{\pd z_2}\Phi_{g-1; n_1+1, n_2, n_3}(z_1, z_2,x^1_{[n_1]_i}; x^2_{[n_2]}; x^3_{[n_3]}; a)|_{z_1, z_2=x_k^1}\\
& = & \half (-a(a+1))^{n_1+n_2+n_3} \sum_{k=1}^{n_1} \sum_{b,c,b^i_j \geq 0}
\cor{\tau_b \tau_c \prod_{\substack{1 \leq i \leq 3, 1 \leq j \leq n_i \\ (i,j) \neq (1,k)}}
\tau_{b^i_j} \cdot T_g(a)}_g \\
&& \cdot \psi_{b+1}(v^1_k; a) \psi_{c+1}(v^1_k;a) \cdot
\prod_{\substack{1 \leq i \leq 3, 1 \leq j \leq n_i \\ (i,j) \neq (1,k)}} \xi^i_{b^i_j}(v^i_j; a_i).
\een
The second term on the right-hand side has several cases.
Case 1. The splitting is stable, i.e.,
\ben
2g_1-1 + |A^1|+|A^2|+|A^3| & > 0, \\
2g_2-1+|B^1|+|B^2|+|B^3| & > 0.
\een
Then we get a term of the form:
\ben
&& \frac{1}{2}
x^1_k\frac{\pd}{\pd x^1_k} \Phi_{g_1; |A^1|+1, |A^2|, |A^3|}(x^1_k, x^1_{A^1};x^2_{A^2}; x^3_{A^3}; a) \\
&& \cdot
x^1_k \frac{\pd}{\pd x^1_k} \Phi_{g_2; |B^1|+1, |B^2|, |B^3|}(x^1_k, x^1_{B^1}; x^2_{B^2}; x^3_{B^3}; a) \\
& = & \frac{1}{2}
(-a(a+1))^{n_1+n_2+n_3-1}  \\
&& \cdot \sum_{b,c,b^i_j \geq 0}
\cor{\tau_b\prod_{i=1}^3 \prod_{j \in A^i} \tau_{b^i_j} \cdot T_{g_1}(a)}_{g_1} \cdot
\cor{\tau_c\prod_{i=1}^3 \prod_{j \in B^i} \tau_{b^i_j} \cdot T_{g_2}(a)}_{g_2} \\
&& \cdot \psi_{b+1}(v^1_k; a) \psi_{c+1}(v^1_k;a) \cdot \prod_{\substack{1 \leq i \leq 3, 1 \leq j \leq n_i \\ (i,j) \neq (1,k)}} \xi^i_{b^i_j}(v^i_j; a_i).
\een
Case 2. There are terms that involve exceptional terms of the form:
\ben
&& \sum_{k=1}^{n_1}
x^1_k \frac{\pd}{\pd x^1_k} \Phi_{0; 1, 0, 0}(x^1_k; a) \cdot
x^1_k \frac{\pd}{\pd x^1_k} \Phi_{g; \bn}(x^1_{[n_1]}; x^2_{[n_2]}; x^3_{[n_3]}; a) \\
& = & \sum_{k=1}^{n_1} \xi_{-1}(v^1_k;a) \cdot  (-a(a+1))^{n_1+n_2+n_3-1} \sum_{b^i_j \geq 0}
\cor{\prod_{i=1}^3 \prod_{j=1}^{n_i} \tau_{b^i_j} \cdot T_g(a)}_g \\
&& \cdot \prod_{i=1}^3 \prod_{j=1}^{n_i} \xi^i_{b^i_j+\delta_{i,1}\delta_{j,k}}(v^i_j; a_i).
\een
Here we have used (\ref{eqn:Phi01(1)}).

Case 3.  We have some unstable terms which combined with the third term on the right-hand side gives us:
\ben
&& \sum_{k=1}^{n_1} \sum_{l \in [n_1]_k}
(x^1_k\frac{\pd}{\pd x^1_k} \Phi_{0; 2, 0, 0}(x^1_k, x^1_l; a) - \frac{x^1_j}{x^1_k-x^1_l}) \\
&&  \cdot
x^1_k \frac{\pd}{\pd x^1_k} \Phi_{g; n_1-1, n_2, n_3}(x^1_{[n_1]_l}; x^2_{[n_2]}; x^3_{[n_3]};  a) \\
& =& - \sum_{k=1}^{n_1} \sum_{l \in [n_1]_k} (\frac{1}{y(x^1_k;a)-y(x^1_l;a)} - \frac{1}{y(x^1_k;a)-1} )
x^1_k \frac{\pd y(x^1_k;a)}{\pd x^1_k} \\
&& \cdot (-a(a+1))^{n_1+n_2+n_3-2} \sum_{b^i_j \geq 0}
\cor{\prod_{\substack{1 \leq i \leq 3, 1 \leq j \leq n_i \\ (i,j)\neq (1, l)}} \tau_{b^i_j} \cdot T_g(a)}_g \\
&& \cdot \prod_{\substack{1 \leq i \leq 3, 1 \leq j \leq n_i \\ (i,j)\neq (1, l)}}
\xi^i_{b^i_j+\delta_{i,1}\delta_{j,k}}(v^i_j; a_i).
\een
Here we have used (\ref{eqn:Phi02(1)Der1}).

Case 4.

\ben
&& \sum_{k=1}^{n_1} \sum_{l = 1}^{n_2}
x^1_k \frac{\pd}{\pd x^1_k} \Phi^{0}_{1, 1, 0}(x^1_k; x^2_l; a) \cdot
x^1_k \frac{\pd}{\pd x^1_k} \Phi^{g}_{n_1, n_2-1, n_3}(x^1_{[n_1]}; x^2_{[n_2]_l}; x^3_{[n_3]};  a) \\
& = &  \sum_{k=1}^{n_1} \sum_{l = 1}^{n_2}
(- \frac{x_k^1 \frac{\pd \tilde{y}(x_k^1;a)}{\pd x_k^1}   }{\tilde{y}(x^1_k; a) - y(x_l^2;a_2)}
+ \frac{x_k^1 \frac{\pd \tilde{y}(x_k^1;a)}{\pd x_k^1}}{\tilde{y}(x_k^1;a)-1}) \\
&& \cdot (-a(a+1))^{n_1+n_2+n_3-2} \sum_{b^i_j \geq 0}
\cor{\prod_{\substack{1 \leq i \leq 3, 1 \leq j \leq n_i \\ (i,j)\neq (2, l)}} \tau_{b^i_j} \cdot T_g(a)}_g \\
&& \cdot \prod_{\substack{1 \leq i \leq 3, 1 \leq j \leq n_i \\ (i,j)\neq (2, l)}}
\xi^i_{b^i_j+\delta_{i,1}\delta_{j,k}}(v^i_j; a_i).
\een
Here we have used (\ref{eqn:Phi(12)Der1}).

Case 5.
\ben
&& \sum_{k=1}^{n_1} \sum_{l =1}^{n_3}
x^1_k \frac{\pd}{\pd x^1_k} \Phi^{0}_{1, 0, 1}(x^1_k; x^3_l; a) \cdot
x^1_k \frac{\pd}{\pd x^1_k} \Phi^{g}_{n_1, n_2, n_3-1}(x^1_{[n_1]}; x^2_{[n_2]}; x^3_{[n_3]_l};  a) \\
& = &  \sum_{k=1}^{n_1} \sum_{l = 1}^{n_2}
 \frac{x^1_k \frac{\pd y(x^1_k ;a)}{\pd x^1_k}}{\tilde{y}(x^3_l; a_3) - y(x_k^1;a_1)}
 \cdot (-a(a+1))^{n_1+n_2+n_3-2} \\
 && \cdot \sum_{b^i_j \geq 0}
\cor{\prod_{\substack{1 \leq i \leq 3, 1 \leq j \leq n_i \\ (i,j)\neq (3, l)}} \tau_{b^i_j} \cdot T_g(a)}_g
\cdot \prod_{\substack{1 \leq i \leq 3, 1 \leq j \leq n_i \\ (i,j)\neq (3, l)}}
\xi^i_{b^i_j+\delta_{i,1}\delta_{j,k}}(v^i_j; a_i).
\een
Here we have used (\ref{eqn:Phi(31)Der1}).

So far we have only considered the terms corresponding to cut-and-join in $x^1_k$ variables.
The terms for $x^2_k$ and $x^3_k$ can be obtained by cyclic symmetry.

As in \cite{Zho},
we regard both sides of the equation (\ref{eqn:CJinX}) as meromorphic functions in $v_1^1$,
take the principal parts and then take only the even powers in $v_1^1$.
The left-hand side and the $\dots$ terms have no contributions.
So we get:
\ben
&& \xi^o_{-1}(v^1_1;a) \cdot   \sum_{b^i_j \geq 0}
\cor{\prod_{i=1}^3 \prod_{j=1}^{n_i} \tau_{b^i_j} \cdot T_g(a)}_g
\cdot \prod_{i=1}^3 \prod_{j=1}^{n_i} \xi^i_{b^i_j+\delta_{i,1}\delta_{j,1}}(v^i_j; a_i) \\
& = &  \half a(a+1) \sum_{b,c,b^i_j \geq 0}
\cor{\tau_b \tau_c \prod_{\substack{1 \leq i \leq 3, 1 \leq j \leq n_i \\ (i,j) \neq (1,1)}}
\tau_{b^i_j} \cdot T_g(a)}_g \\
&& \cdot \psi_{b+1}(v^1_1; a) \psi_{c+1}(v^1_1;a) \cdot
\prod_{\substack{1 \leq i \leq 3, 1 \leq j \leq n_i \\ (i,j) \neq (1,1)}} \xi^i_{b^i_j}(v^i_j; a_i) \\
& - &  \frac{1}{2} \sum^{stable} \sum_{b,c,b^i_j \geq 0}
\cor{\tau_b\prod_{i=1}^3 \prod_{j \in A^i} \tau_{b^i_j} \cdot T_{g_1}(a)}_{g_1} \cdot
\cor{\tau_c\prod_{i=1}^3 \prod_{j \in B^i} \tau_{b^i_j} \cdot T_{g_2}(a)}_{g_2} \\
&& \cdot \psi_{b+1}(v^1_1; a) \psi_{c+1}(v^1_1;a) \cdot
\prod_{\substack{1 \leq i \leq 3, 1 \leq j \leq n_i \\ (i,j) \neq (1,1)}} \xi^i_{b^i_j}(v^i_j; a_i) \\
& - & \frac{1}{a(a+1)} \sum_{l \in [n_1]_1} (\frac{1}{y(x^1_1;a)-y(x^1_l;a)} - \frac{1}{y(x^1_1;a)-1} )
x^1_1 \frac{\pd y(x^1_1;a)}{\pd x^1_1} \\
&& \cdot  \sum_{b^i_j \geq 0}
\cor{\prod_{\substack{1 \leq i \leq 3, 1 \leq j \leq n_i \\ (i,j)\neq (1, l)}} \tau_{b^i_j} \cdot T_g(a)}_g
\cdot \prod_{\substack{1 \leq i \leq 3, 1 \leq j \leq n_i \\ (i,j)\neq (1, l)}}
\xi^i_{b^i_j+\delta_{i,1}\delta_{j,1}}(v^i_j; a_i) \\
& + &  \frac{1}{a(a+1)} \sum_{l = 1}^{n_2}
(- \frac{x_1^1 \frac{\pd \tilde{y}(x_1^1;a)}{\pd x_1^1}   }{\tilde{y}(x^1_1; a) - y(x_l^2;a_2)}
+ \frac{x_1^1 \frac{\pd \tilde{y}(x_1^1;a)}{\pd x_1^1}}{\tilde{y}(x_1^1;a)-1}) \\
&& \cdot  \sum_{b^i_j \geq 0}
\cor{\prod_{\substack{1 \leq i \leq 3, 1 \leq j \leq n_i \\ (i,j)\neq (2, l)}} \tau_{b^i_j} \cdot T_g(a)}_g
\cdot \prod_{\substack{1 \leq i \leq 3, 1 \leq j \leq n_i \\ (i,j)\neq (2, l)}}
\xi^i_{b^i_j+\delta_{i,1}\delta_{j,1}}(v^i_j; a_i) \\
& + & \frac{1}{a(a+1)} \sum_{l = 1}^{n_2}
 \frac{x^1_1 \frac{\pd y(x^1_1 ;a)}{\pd x^1_k}}{\tilde{y}(x^3_l; a_3) - y(x_1^1;a_1)} \\
 && \cdot \sum_{b^i_j \geq 0}
\cor{\prod_{\substack{1 \leq i \leq 3, 1 \leq j \leq n_i \\ (i,j)\neq (3, l)}} \tau_{b^i_j} \cdot T_g(a)}_g
\cdot \prod_{\substack{1 \leq i \leq 3, 1 \leq j \leq n_i \\ (i,j)\neq (3, l)}}
\xi^i_{b^i_j+\delta_{i,1}\delta_{j,1}}(v^i_j; a_i),
\een
modulo terms analytic in $v_1$.
We now take $\prod_{j=2}^{n_1} x^1_j \frac{\pd}{\pd x^1_j} \prod_{i=2}^3 \prod_{j=1}^{n_i} x^i_j\frac{\pd}{\pd x^i_j}$
on both sides then dividing both sides by $\xi_{-1}^o(v;a)$.
This gives us:
\ben
&&  \sum_{b^i_j \geq 0}
\cor{\prod_{i=1}^3 \prod_{j=1}^{n_i} \tau_{b^i_j} \cdot T_g(a)}_g
\cdot \prod_{i=1}^3 \prod_{j=1}^{n_i} \xi^i_{b^i_j+1}(v^i_j; a_i) \\
& = & \frac{1}{\xi^o_{-1}(v^1_1;a)} \cdot \biggl(   \half a(a+1) \sum_{b,c,b^i_j \geq 0}
\cor{\tau_b \tau_c \prod_{\substack{1 \leq i \leq 3, 1 \leq j \leq n_i \\ (i,j) \neq (1,1)}}
\tau_{b^i_j} \cdot T_g(a)}_g \\
&& \cdot \psi_{b+1}(v^1_1; a) \psi_{c+1}(v^1_1;a) \cdot
\prod_{\substack{1 \leq i \leq 3, 1 \leq j \leq n_i \\ (i,j) \neq (1,1)}} \xi^i_{b^i_j + 1}(v^i_j; a_i) \\
& - &  \frac{1}{2} \sum^{stable} \sum_{b,c,b^i_j \geq 0}
\cor{\tau_b\prod_{i=1}^3 \prod_{j \in A^i} \tau_{b^i_j} \cdot T_{g_1}(a)}_{g_1} \cdot
\cor{\tau_c\prod_{i=1}^3 \prod_{j \in B^i} \tau_{b^i_j} \cdot T_{g_2}(a)}_{g_2} \\
&& \cdot \psi_{b+1}(v^1_1; a) \psi_{c+1}(v^1_1;a) \cdot
\prod_{\substack{1 \leq i \leq 3, 1 \leq j \leq n_i \\ (i,j) \neq (1,1)}} \xi^i_{b^i_j +1}(v^i_j; a_i) \\
& - & \frac{1}{a(a+1)} \sum_{l \in [n_1]_1} \frac{1}{(y(x^1_1;a)-y(x^1_l;a))^2}
\cdot x^1_1 \frac{\pd y(x^1_1;a)}{\pd x^1_1} \cdot x^1_l \frac{\pd y(x^1_l;a)}{\pd x^1_l}  \\
&& \cdot  \sum_{b^i_j \geq 0}
\cor{\prod_{\substack{1 \leq i \leq 3, 1 \leq j \leq n_i \\ (i,j)\neq (1, l)}} \tau_{b^i_j} \cdot T_g(a)}_g
\cdot \prod_{\substack{1 \leq i \leq 3, 1 \leq j \leq n_i \\ (i,j)\neq (1, l)}}
\xi^i_{b^i_j+1}(v^i_j; a_i) \\
& - &  \frac{1}{a(a+1)} \sum_{l = 1}^{n_2} \frac{1}{(\tilde{y}(x^1_1; a) - y(x_l^2;a_2))^2}
\cdot x_1^1 \frac{\pd \tilde{y}(x_1^1;a)}{\pd x_1^1}  \cdot x_l^2 \frac{\pd y(x_l^2;a_2)}{\pd x_l^2}   \\
&& \cdot  \sum_{b^i_j \geq 0}
\cor{\prod_{\substack{1 \leq i \leq 3, 1 \leq j \leq n_i \\ (i,j)\neq (2, l)}} \tau_{b^i_j} \cdot T_g(a)}_g
\cdot \prod_{\substack{1 \leq i \leq 3, 1 \leq j \leq n_i \\ (i,j)\neq (2, l)}}
\xi^i_{b^i_j+1}(v^i_j; a_i) \\
& - & \frac{1}{a(a+1)} \sum_{l = 1}^{n_2}\frac{1}{(y(x_1^1;a_1)-\tilde{y}(x^3_l; a_3))^2}
\cdot x^1_1 \frac{\pd y(x^1_1 ;a)}{\pd x^1_k} \cdot x_l^3 \frac{\pd \tilde{y}(x_l^3;a)}{\pd x_l^3}  \\
&& \cdot \sum_{b^i_j \geq 0}
\cor{\prod_{\substack{1 \leq i \leq 3, 1 \leq j \leq n_i \\ (i,j)\neq (3, l)}} \tau_{b^i_j} \cdot T_g(a)}_g
\cdot \prod_{\substack{1 \leq i \leq 3, 1 \leq j \leq n_i \\ (i,j)\neq (3, l)}}
\xi^i_{b^i_j+1}(v^i_j; a_i) \biggr),
\een
modulo  a term with at most a simple pole at $0$ in $v_1^1$.

\section{Eynard-Orantin Recursion relations for Three-Partition Triple Hodge Integrals}

In this section we  reformulate the recursion relations for three-partition triple Hodge integrals
derived in the end of last section as Eynard-Orantin type recursion relations.
This verifies a version of local mirror symmetry proposed by Bouchard-Klemm-Mari\~no-Pasquetti
\cite{BKMP} for the topological vertex.

\subsection{Differentials associated to  three-partition triple Hodge integrals}

Define
\begin{multline}
W_g(x^1_{[n_1]}; x^2_{[n_2]}; x^3_{[n_3]};a) \\
= (-1)^{g-1} \prod_{i=1}^3 \prod_{j=1}^{n_i} \frac{\pd}{\pd x^i_j}
\Phi^g_{n_1,n_2, n_3}(x^1_{[n_1]}; x^2_{[n_2]}; x^3_{[n_3]}; a)
\cdot \prod_{i=1}^3 \prod_{j=1}^{n_i} dx^i_j.
\end{multline}
Then we have
\begin{multline} \label{eqn:Wg}
W_g(x^1_{[n_1]}; x^2_{[n_2]}; x^3_{[n_3]}; a)
 = (-1)^{g+n_1+n_2+n_3} (a(a+1))^{n_1+n_2+n_3-1} \\
\cdot \cor{\prod_{i=1}^3 \prod_{j=1}^{n_i} \tau_{b^i_j} \cdot T_g(a)}_g
\prod_{i=1}^3 \prod_{j=1}^{n_i} d \phi_{b^i_j}(x^i_j; a_i)
\cdot \prod_{j=1}^{n_2}  \frac{1}{a^{2+b^2_j}} \cdot \prod_{j=1}^{n_3}  \frac{1}{(-a-1)^{2+b^3_j}},
\end{multline}
for $2g-2+n > 0$.
By (\ref{eqn:Phi01(1)}), (\ref{eqn:Phi01(2)}) and (\ref{eqn:Phi01(3)}),
we have
\bea
&& W_0(x_1^1;a) = \ln y(x_1^1;a) {d x_1^1 \over x^1_1}, \label{eqn:W0y(1)} \\
&& W_0(x_1^2;a) = \ln y(x_1^2;a_2) {d x_1^2 \over x^2_1}, \label{eqn:W0y(2)} \\
&& W_0(x_1^3;a) = \ln y(x_1^3;a_3) {d x_1^3 \over x^3_1}. \label{eqn:W0y(3)}
\eea
By (\ref{eqn:Phi02(1)}) we get:
\be
W_0(x^1_1,x^1_2;a) = \frac{dy(x^1_1;a) dy(x^1_2;a)}{(y(x^1_1;a)-y(x^1_2;a))^2}
- \frac{dx^1_1 dx^1_2}{(x^1_1-x^1_2)^2}. \label{eqn:W0yy(1)}
\ee
Similarly,
\bea
&& W_0(x^2_1,x^2_2;a) = \frac{dy(x^2_1;a_2) dy(x^2_2;a_2)}{(y(x^2_1;a_2)-y(x^2_2;a_2))^2}
- \frac{dx^2_1 dx^2_2}{(x^2_1-x^2_2)^2}, \label{eqn:W0yy(2)} \\
&& W_0(x^3_1,x^3_2;a) = \frac{dy(x^3_1;a_3) dy(x^3_2;a_3)}{(y(x^3_1;a_3)-y(x^3_2;a_3))^2}
- \frac{dx^3_1 dx^3_2}{(x^3_1-x^3_2)^2}. \label{eqn:W0yy(3)}
\eea
By (\ref{eqn:Phi(12)Der2}), (\ref{eqn:Phi(23)Der2}) and (\ref{eqn:Phi(31)Der2}),
we have
\bea
&& W_0(x_1^1; x_1^2; a)
= \frac{d \tilde{y}(x^1_1;a) d y(x_1^2;a_2) }{(\tilde{y}(x^1_1;a) - y(x_1^2;a_2))^2}, \label{eqn:W0yy(12)} \\
&& W_0(x_1^2; x_1^3; a)
= \frac{d \tilde{y}(x^2_1;a_2) d y(x_1^3;a_3)}{(\tilde{y}(x^2_1;a_2) - y(x_1^3;a_3))^2}, \label{eqn:W0yy(23)} \\
&& W_0(x_1^1; x_1^3; a)
= \frac{d \tilde{y}(x^3_1;a_3) d y(x_1^1;a)}{(\tilde{y}(x^3_1;a_3) - y(x_1^1;a))^2}. \label{eqn:W0yy(31)}
\eea

\subsection{Eynard-Orantin formalism for the topological vertex}

According to the proposal in \cite{Mar, BKMP},
the differentials $W_g(x^1_{[n_1]}; x^2_{[n_2]}; x^3_{[n_3]};a)$ can be computed recursively
by the Enyard-Orantin formalism,
for the framed mirror curve given by (\ref{eqn:FramedMirror}).
By abuse of notations,
we will write the differential $W_g(x^1_{[n_1]}; x^2_{[n_2]}; x^3_{[n_3]};a)$ as $W_g(y^1_{[n_1]}; y^2_{[n_2]}; y^3_{[n_3]};a)$,
where $y^i_j = y(x^i_j;a_i)$.
The initial values are given by (\ref{eqn:W0y(1)})-(\ref{eqn:W0yy(31)}),
and the recursion is given by:
\begin{multline} \label{eqn:EO}
W_g(y^1_{[n_1]}; y^2_{[n_2]}; y^3_{[n_3]}; a)
= \res_{z=0}{d E_{z}(y^1_1) \over \omega(z)} \\
\quad  \Big( W_{g-1} (y(z), y(P(z)), y^1_{[n_1]_1}; y^2_{[n_2]}; y^3_{[n_3]}; a ) \\
 \quad +\sum_{\substack{g_1+g_2=g\\ A^1\coprod B^1 = [n]_1}}
 \sum_{\substack{A^2 \coprod B^2 =[n_2] \\A^3 \coprod B^3 = [n_3]}}
 W_{g_1}(y(z), y^1_{A^1}; y^2_{A^2}; y^3_{A^3}; a) \\
\cdot W_{g_2}(y(P(z)), y^1_{B^1}; y^2_{B^2}; y^3_{B^3}; a)\Big),
\end{multline}
where
\bea
&& \omega(z) =  ( \ln y(z) - \ln y(P(z)) )  \cdot {d x(z) \over x(z) }, \\
&& d E_z (y^1_1) = {1 \over 2}\left( {1 \over y(z^1_1) - y(z)} - {1 \over y(z^1_1) - y(P(z))} \right) dy_1.
\eea
Here for simplicity of notations,
we write
\bea
&& x(z) = x(\frac{a}{a+1}+z;a) = (\frac{1}{a+1} - z) (\frac{a}{a+1} + z)^a, \\
&& y(z) = \frac{a}{a+1} + z, \qquad y^1_1 = y(z^1_1).
\eea
On the right-hand side of (\ref{eqn:EO}),
$W_0(\frac{a}{a+1}+z;a)$ and $W_0(\frac{a}{a+1}+P(z);a)$ are understood as $0$.
Therefore,
(\ref{eqn:EO}) is a recursion relation that determines
$W_g(y^1_{[n_1]}; y^2_{[n_2]}; y^3_{[n_3]};a)$ for $n_1 > 0$.
When $n_1 = 0$,
one can use the cyclic symmetry (\ref{eqn:Cyclic})  to reduce to the $n_1 > 0$, $n_3 = 0$  case.

Theorem 1 in the Introduction can be rephrased more precisely as

\begin{theorem}
The differentials $W_g(y^1_{[n_1]}; y^2_{[n_2]}; y^3_{[n_3]}; a)$ satisfy the recursion relations
(\ref{eqn:EO}).
\end{theorem}

When $n_2 =n_3 = 0$,
this has been proved in \cite{Zho} and \cite{Che} using ideas
from \cite{EMS}.
For the proof of the general case,
one needs to show that (\ref{eqn:EO}) is equivalent to the recursion relations derived from the cut-and-join equation
in the end of \S \ref{sec:CJ}.
It is almost a verbatim  straightforward generalization of the treatment in \cite{Zho},
so we will omit it.

\end{document}